\documentclass[12pt]{amsart}

\usepackage{array}

\usepackage{graphicx}
\usepackage{float}

\usepackage{amssymb,amsmath}
\usepackage{amssymb,latexsym}
\usepackage{amsfonts}
\usepackage{amssymb}
\usepackage{stmaryrd}
\usepackage{longtable}
\usepackage{amsmath, geometry, amssymb} 
\usepackage{chngcntr}
\numberwithin{equation}{section}
\newtheorem{theorem}{Theorem}[section]
\newtheorem{lemma}{Lemma}[section]
\newtheorem{remark}{Remark}[section]
\newtheorem{corollary}{Corollary}[section]

\newcommand{\pmdd}{\hspace{-3mm} \pmod}

\newcommand{\sqr}[2]{{\vcenter{\vbox{\hrule height#2pt
                \hbox{\vrule width#2pt height#1pt \kern#1pt
                \vrule width#2pt}\hrule height#2pt}}}}

\newcommand{\mb}{\mbox}

\newcommand{\beq}{\begin{equation}}
\newcommand{\eeq}{\end{equation}}
\newcommand{\beqar}{\begin{eqnarray}}
\newcommand{\eeqar}{\end{eqnarray}}
\def\beqars{\begin{eqnarray*}}
\def\eeqars{\end{eqnarray*}}
\def\eop{\hfill\mb{$\hspace{12pt}\vrule height 7pt width 6pt depth 0pt$}}

\newcommand{\pmd}{\hspace{-3mm} \pmod}

\newcommand{\smod}[1]{\hspace{-1mm} \pmod{#1}}

\newcommand{\qu}[2]{\Bigl({\frac{#1}{#2}}\Bigr) }
\newcommand{\dqu}[2]{\ds{\qu{#1}{#2}}}

\def \ds{\displaystyle}

\newcommand{\nn}{\mathbb{N}}
\newcommand{\zz}{\mathbb{Z}}

\newcommand{\cc}{\mathbb{C}}

\newcommand{\hh}{\mathbb{H}}

\allowdisplaybreaks

\begin{document}

\title{Representations by Quaternary Quadratic Forms with Coefficients  $1$, $2$, $5$ or $10$}

\author{Ay\c{s}e Alaca and Mada Altiary}

\maketitle

\markboth{Ay\c{s}e Alaca and Mada Altiary}
{Representations by Quaternary Quadratic Forms with Coefficients  $1$, $2$, $5$ or $10$}

\begin{abstract}
We determine explicit formulas for the number of representations of a positive integer $n$ by quaternary quadratic forms 
with coefficients $1$, $2$, $5$ or $10$. We use  a modular forms approach.

\vspace{2mm}


\vspace{2mm}

\end{abstract}

\section{Introduction}

Let $\nn$, $\nn_0$, $\zz$ and $\cc$ denote the sets of
positive integers, nonnegative integers, integers and complex numbers, respectively. For $n\in\nn$ we set
$\sigma(n)= \sum_{ d \mid n } d$, where $d$ runs through the positive divisors of $n$.
If $n\not \in\nn$ we set $\sigma(n)=0$.
For $a_1, a_2, a_3, a_4 \in \nn$,  and $n \in \nn_0$ 
we define
\begin{align*}
N(a_1,a_2,a_3,a_4 ;n):={\rm card}\{(x_1,x_2,x_3,x_4)\in \zz^4 \mid n= a_1x_1^2 +a_2x_2^2+a_3x_3^2+ a_4 x_4^2 \}.
\end{align*}
It is a classical result of Jacobi  \cite {jacobi, williamsBook}  that
\begin{align*}
N(1,1,1,1;n)=8\sigma(n) - 32 \sigma (n/4).
\end{align*}
Formulas for $N(a_1, a_2, a_3, a_4;n)$ for the quaternary quadratic forms
\begin{align*}
(a_1,a_2,a_3,a_4)= (1,1,1,2), (1,1,2,2), (1,2,2,2), (1,1,1,5), (1,1,5,5), (1,5,5,5)
\end{align*}
are in the literature, 
see for example \cite{alaca2007quaternary,  alaca2007nineteen, ASLW-2009,  K1155, L1122, L1112-1222, L1115, L1555, L1155, Lomadse, W1112-1222}.

There are twenty-six quaternary quadratic forms 
$ a_1x_1^2 +a_2x_2^2+a_3x_3^2+ a_4 x_4^2 $, where $a_1, a_2, a_3, a_4 \in \{1,2,5,10\}$, $a_1\leq a_2\leq a_3\leq a_4$ and 
$\gcd (a_1,a_2,a_3,a_4)=1$, see Table 2.1.  
In this paper, we  determine  an explicit formula for $N(a_1,a_2,a_3,a_4 ;n)$ for each of these quaternary 
forms in a uniform manner. We use  a modular forms approach. 

For $q \in \cc$ with $|q|<1$, Ramanujan's theta function $\varphi (q)$  is defined by 
\begin{align*}
 \varphi (q) = \sum_{n=-\infty}^{\infty} q^{n^2}.
\end{align*}
For $a_1, a_2,a_3,a_4 \in \nn$ we have
\begin{align}
\sum_{n=1}^{\infty} N(a_1,a_2,a_3,a_4;n) q^n=\varphi(q^{a_1}) \varphi(q^{a_2})\varphi(q^{a_3})\varphi(q^{a_4}). 
\end{align}
The Dedekind eta function $\eta (z)$ is the holomorphic function defined on the upper half plane $\hh = \{ z \in \cc \mid \mbox{\rm Im}(z) >0 \}$ 
by 
\begin{align*}
\eta (z) = e^{\pi i z/12} \prod_{n=1}^{\infty} (1-e^{2\pi inz}).
\end{align*}
Throughout the remainder of  the paper we take $q=q(z):=e^{2\pi i z}$ with $z\in \hh$. 
Thus we  can express $\eta (z)$ as 
\begin{align}
\eta (z) = q^{1/24} \prod_{n=1}^{\infty} (1-q^{n}).
\end{align} 
An eta quotient is defined to be a finite product of the form
\begin{align*}
f(z) = \prod_{\delta } \eta^{r_{\delta}} ( \delta z),
\end{align*}
where $\delta$ runs through a finite set of positive integers and the exponents $r_{\delta}$ are non-zero integers.
It is known (see for example \cite[Corollary 1.3.4]{Berndt}) that  
\begin{align}
\varphi(q) = \frac{\eta^5(2z)}{\eta^2(z) \eta^2(4z)}. 
\end{align}

\section{Modular spaces $M_2(\Gamma_0(40),\chi_i)$ with $i\in\{0,1,2,3\}$}  

For $n\in \nn$ and Dirichlet characters $\chi$ and $\psi$  we define  $\displaystyle \sigma_{\chi,\psi}(n)$ by
\begin{align}
\sigma_{\chi,\psi }(n) :=\sum_{1 \leq m|n}\psi(m)\chi(n/m)m.
\end{align}
If $n \not\in \nn$ we set $\sigma_{\chi,\psi }(n)=0$. Let $\chi_0$  denote the trivial character, that is $\chi_0 (n) =1$ for all $n \in \zz$. 
Hence  $\sigma_{\chi_0, \chi_0}(n) $ coincides with  the sum of divisors function $\sigma (n)$. Let $N\in\nn$.  The modular  subgroup $\Gamma_0(N)$ is defined by
\begin{align*}
\Gamma_0(N) = \left\{ \left(
\begin{array}{cc}
a & b \\
c & d
\end{array}
\right) \Big | \;  a,b,c,d\in \zz ,~ ad-bc = 1,~c \equiv 0 \pmd {N}
\right\} .
\end{align*}
Let $\chi$ be a Dirichlet character of modulus dividing $N$ and let $k\in \zz$.  
We write $M_k(\Gamma_0(N),\chi)$ to denote the space of modular forms of weight $k$ with multiplier system $\chi$ for $\Gamma_0(N)$, and $E_{k}(\Gamma_0(N),\chi)$ and $S_{k}(\Gamma_0(N),\chi)$ to denote the subspaces of Eisenstein forms and cusp forms of $M_{k}(\Gamma_0(N),\chi)$, respectively. If $\chi =\chi_0$, then we write $M_k(\Gamma_0(N))$  for $ M_k(\Gamma_0(N),\chi_0)$, and $S_k(\Gamma_0(N))$ for $S_k(\Gamma_0(N),\chi_0)$. 
It is known (see for example \cite[p. 83]{stein}) that
\begin{align}
M_{k}(\Gamma_0(N),\chi)=E_{k}(\Gamma_0(N),\chi)\oplus S_{k}(\Gamma_0(N),\chi).
\end{align}
For $n \in \zz$ we define three  Dirichlet characters  by 
\begin{align}
\chi_1 ( n)=\dqu{5}{n},~\chi_2 (n) =\dqu{8}{n},~\chi_3 (n) =\dqu{40}{n}.
\end{align}
We define the Eisenstein series 
\begin{align}
& L(q):=E_{\chi_0,\chi_0}(q)=-\frac{1}{24}+ \sum_{n=1}^{\infty} \sigma(n) q^n,\\
& E_{\chi_0,\chi_1}(q)=-\frac{1}{5}+\sum_{n=1}^{\infty} \sigma_{\chi_0, \chi_1}(n) q^n,~
 E_{\chi_1,\chi_0}(q)=\sum_{n=1}^{\infty} \sigma_{\chi_1, \chi_0}(n) q^n,\\
&E_{\chi_0,\chi_2}(q)=-\frac{1}{2}+\sum_{n=1}^{\infty} \sigma_{\chi_0, \chi_2}(n) q^n,~
 E_{\chi_2,\chi_0}(q)=\sum_{n=1}^{\infty} \sigma_{\chi_2, \chi_0}(n) q^n,\\
&E_{\chi_0,\chi_3}(q)=-7+\sum_{n=1}^{\infty} \sigma_{\chi_0, \chi_3}(n) q^n, ~
 E_{\chi_3,\chi_0}(q)=\sum_{n=1}^{\infty} \sigma_{\chi_3, \chi_0}(n) q^n,\\
& E_{\chi_1,\chi_2}(q)=\sum_{n=1}^{\infty} \sigma_{\chi_1, \chi_2}(n) q^n,~\hspace{10mm}
 E_{\chi_2,\chi_1}(q)=\sum_{n=1}^{\infty} \sigma_{\chi_2, \chi_1}(n) q^n.
\end{align}

We use the following lemma to determine if certain eta quotients are modular forms.
See \cite[p. 174]{Gordon},  \cite[Corollary 2.3, p. 37]{Kohler}, \cite[Theorem 5.7, p. 101]{Kilford}, \cite{Ligozat} and \cite[Theorem 1.64]{ono}.
 
\begin{lemma} {\rm \bf  (Ligozat)}  
Let $N\in\nn$ and  $f(z)=\ds \prod_{1\leq \delta \mid N} \eta^{r_{\delta}}(\delta z)$ be an eta quotient.
Let $s=\ds \prod_{ 1\leq \delta \mid N} \delta^{|r_{\delta}|}$  and  $\ds k = \frac{1}{2} \sum_{1 \leq \delta \mid N} r_{\delta}$.
Suppose that the following conditions are satisfied:\\

{\em (L1)~} $\ds \sum_{ 1\leq  \delta \mid N} \delta \cdot r_{\delta} \equiv 0 \smod {24}$,\\

{\em (L2)~} $\ds \sum_{ 1 \leq \delta \mid N} \frac{N}{\delta} \cdot r_{\delta} \equiv 0 \smod {24}$,\\

{\em (L3)~} $\ds \sum_{1 \leq \delta \mid N} \frac{ \gcd (d, \delta)^2 \cdot r_{\delta} }{\delta} \geq 0 $  for each 
positive divisor $d$ of $N$, \\

{\em (L4)~}  $k$ is an  integer.

\noindent
Then $f(z) \in M_k(\Gamma_0(N),\chi)$, where the character $\chi$ is given by  $\chi (m) = \ds \dqu{(-1)^ks}{m}$.

{\em (L3)$'$~} In addition to the above conditions, if the inequality in {\em(L3)} is strict  for each positive divisor $d$ of $N$,
then $f(z) \in S_k(\Gamma_0(N), \chi)$.
\end{lemma}

In Table 2.1, we group our twenty-six quaternary forms $(a_1, a_2, a_3 ,a_4)$  
according to modular spaces $M_2(\Gamma_0(40),\chi)$ to which 
$\varphi(q^{a_1})\varphi(q^{a_2})\varphi(q^{a_3})\varphi(q^{a_4})$ belong.

\begin{center}
Table 2.1\\[1mm]
\begin{tabular}{|l|l|l|l|} \hline
$M_2(\Gamma_0(40))$ & $ M_2(\Gamma_0(40),\chi_1)$ & $M_2(\Gamma_0(40),\chi_2)$ & $M_2(\Gamma_0(40),\chi_3)$\\
\hline
$(1,1,1,1) \checkmark$  & $(1,1,1,5)\checkmark$      & $(1,1,1,2)\checkmark$     &  $(1,1,1,10)$  \\
$(1,1,2,2)\checkmark$  &  $(1,1,2,10)\ast$              &  $(1,1,5,10)$                   &$(1,1,2,5)\ast$  \\
$(1,1,5,5)\checkmark$  &$(1,2,2,5)\ast$                 & $(1,2,2,2)\checkmark$       & $(1,2,2,10)$  \\
$(1,1,10,10)$              & $(1,5,5,5)\checkmark$     & $(1,2,5,5)$                     & $(1,5,5,10)$ \\
$(1,2,5,10)\ast$          & $(1,5,10,10)$                  & $(1,2,10,10)$                   &  $(1,10,10,10)$   \\
$(2,2,5,5)$               & $(2,5,5,10)$                     & $(2,2,5,10)$                   & $(2,2,2,5)$  \\
                                &                                       &                                        & $(2,5,5,5)$ \\                
                                &                                      &                                        & $(2,5,10,10)$  \\
\hline
\end{tabular}
\end{center}
Formulas $N(a_1,a_2,a_3,a_4;n)$ for the forms  with a checkmark ($\checkmark$) in Table 2.1  are known.
Of the remaining nineteen forms, four are universal and identified with an asterisk ($\ast$).

We deduce from \cite[Sec. 6.1, p. 93]{stein} that
\begin{align}
\dim(E_2(\Gamma_0(40)))=7,~\dim(S_2(\Gamma_0(40)))=3.
\end{align}
We also deduce from \cite[Sec. 6.3, p. 98]{stein} that
\begin{align}
&\dim(E_2(\Gamma_0(40),\chi_1))=8,~\dim(S_2(\Gamma_0(40),\chi_1))=2, \\
&\dim(E_2(\Gamma_0(40),\chi_2))=4,~\dim(S_2(\Gamma_0(40),\chi_2))=4, \\
&\dim(E_2(\Gamma_0(40),\chi_3))=4,~\dim(S_2(\Gamma_0(40),\chi_3))=4. 
\end{align}

\begin{theorem}  
Let $\chi_1, \chi_2, \chi_3$ be as in {\em (2.3)}. 
If $(a_1,a_2,a_3,a_4)$  is in the first, second, third or fourth column of {\em Table 2.1}, then 
\begin{align*} 
&\varphi(q^{a_1})\varphi(q^{a_2})\varphi(q^{a_3})\varphi(q^{a_4}) \in M_2(\Gamma_0(40)),\\
&\varphi(q^{a_1})\varphi(q^{a_2})\varphi(q^{a_3})\varphi(q^{a_4}) \in M_2(\Gamma_0(40),\chi_1),\\
&\varphi(q^{a_1})\varphi(q^{a_2})\varphi(q^{a_3})\varphi(q^{a_4}) \in M_2(\Gamma_0(40),\chi_2),\\
&\varphi(q^{a_1})\varphi(q^{a_2})\varphi(q^{a_3})\varphi(q^{a_4}) \in M_2(\Gamma_0(40),\chi_3),
\end{align*}
respectively.
\end{theorem}

{\bf Proof.} 
The assertion directly follows from (1.3) and  Lemma 2.1.
\eop

\vspace{1mm}

Let $n\in\nn$. We define the eta quotients $A_k(q)$,  $B_k(q)$,  $C_k(q)$, $D_k(q)$
and integers $a_k(n)$, $b_k(n)$,  $c_k(n)$, $d_k(n)$ as follows: 
\begin{align}
&A_1(q)=\displaystyle\sum_{n=1}^{\infty} a_1(n)q^n = \eta^2(2z)\eta^2({10}z),\\
&A_2(q)=A_1(q^2)=\displaystyle\sum_{n=1}^{\infty} a_2(n)q^n = \eta^2(4z)\eta^2({20z}), \\
&A_3(q)=\displaystyle\sum_{n=1}^{\infty} a_3(n)q^n = \frac{\eta^5(4z)\eta({10}z)\eta^2({40}z)}{\eta(2z)\eta^2(8z)\eta({20}z)},\\
&B_1(q)=\sum_{n=1}^{\infty}b_1(n)q^n=\frac{\eta(2z) \eta^4(20 z)}{\eta(10z)},\\
&B_2(q)=\sum_{n=1}^{\infty}b_2(n)q^n=\frac{\eta^4(4z) \eta(10 z)}{\eta(2z)}, \\
&C_1(q)=\sum_{n=1}^{\infty}c_1(n)q^n=\frac{\eta^2(z) \eta(8 z)\eta^2(10 z)\eta(40 z)}{\eta(2z)\eta(20z)},\\
&C_2(q)=\sum_{n=1}^{\infty}c_2(n)q^n=\frac{\eta(z) \eta(5z)\eta^2(8z)\eta^2(20z)}{\eta(4z)\eta(10z)},\\
&C_3(q)=\sum_{n=1}^{\infty}c_3(n)q^n=\frac{\eta^6(2z) \eta(10z)\eta^2(40z)}{\eta^2(z)\eta^2(4z)\eta(20z)},\\
&C_4(q)=\sum_{n=1}^{\infty}c_4(n)q^n=\frac{\eta^6(4z) \eta^2(5z)\eta(20z)}{\eta^2(2z)\eta^2(8z)\eta(10z)},\\
&D_1(q)=\sum_{n=1}^{\infty}d_1(n)q^n=\frac{\eta^2(z) \eta^6(4 z)\eta(20 z)}{\eta^3(2z)\eta^2(8z)},\\
&D_2(q)=\sum_{n=1}^{\infty}d_2(n)q^n=\frac{\eta^2(5z) \eta(8z)\eta(10z)\eta(40z)}{\eta(20z)},\\
&D_3(q)=\sum_{n=1}^{\infty}d_3(n)q^n=\frac{\eta(z) \eta(5z)\eta(20z)\eta^2(40z)}{\eta(10z)},\\
&D_4(q)=\sum_{n=1}^{\infty}d_4(n)q^n=\frac{\eta(z) \eta(4z)\eta(5z)\eta^2(8z)}{\eta(2z)}.
\end{align}

\begin{theorem} 
Let $\chi_1, \chi_2, \chi_3$ be as in {\em (2.3)}. Then 
\begin{align*}
&\{A_1(q), A_2(q), A_3(q)\}, \hspace{12mm}
\{B_1(q), B_2(q)\}, \\
&\{C_1 (q), C_2 (q), C_3 (q), C_4 (q) \}, ~~
\{D_1(q), D_2(q), D_3(q), D_4(q)\}
\end{align*}
are bases for  $S_2(\Gamma_0(40))$, $S_2(\Gamma_0(40),\chi_1)$, $S_2(\Gamma_0(40),\chi_2)$ and $S_2(\Gamma_0(40),\chi_3)$, respectively.
\end{theorem}
 
{\bf Proof.} The set $\{A_1(q), A_2(q), A_3(q)\}$ is linearly independent over $\cc$. 
By Lemma 2.1, we have $A_k(q)\in S_2(\Gamma_0(40))$ for $k=1,2,3$. The assertion now follows from (2.9).
Similarly, the remaining three assertions follow from (2.10), (2.11), (2.12) and  Lemma 2.1.
\eop

\begin{theorem}  
Let $\chi_0$ be the trivial character and $\chi_1, \chi_2, \chi_3$ be as in {\em(2.3)}. Then
\begin{align*}
&\{L(q)-tL(q^t) \mid t=2,4,5,8,10,20,40 \} ,\\
&\{E_{\chi_0,\chi_1}(q^t), E_{\chi_1,\chi_0}(q^t) \mid t=1,2,4,8\}, \\
&\{E_{\chi_0,\chi_2}(q^t), E_{\chi_2,\chi_0}(q^t) \mid t=1,5\}, \\
&\{E_{\chi_0,\chi_3}(q), E_{\chi_1,\chi_2}(q), E_{\chi_2,\chi_1}(q), E_{\chi_3,\chi_0}(q)\} 
\end{align*}
are bases for $E_2(\Gamma_0(40))$,  $ E_2(\Gamma_0(40),\chi_1)$,  $ E_2(\Gamma_0(40),\chi_2)$ and  $ E_2(\Gamma_0(40),\chi_3)$, respectively.
\end{theorem}

{\bf Proof.} 
The assertions follow from \cite[Theorem 5.9]{stein} with
$\chi=\psi=\chi_0$;
$\epsilon=\chi_1$ and $\chi, \psi \in \{\chi_0, \chi_1\}$;
$\epsilon=\chi_2$ and $\chi, \psi \in \{\chi_0, \chi_2\}$;
$\epsilon=\chi_3$ and $\chi, \psi \in \{\chi_0,  \chi_1, \chi_2,  \chi_3\}$, respectively.
\eop

\begin{theorem}  
Let $\chi_0$ be the trivial character and $\chi_1, \chi_2, \chi_3$ be as in {\em(2.3)}. Then 
\begin{align*}
&\{L(q)-tL(q^t)\mid t=2,4,5,8,10,20,40\} \cup \{A_1(q), A_2(q), A_3(q)\}, \\
&\{E_{\chi_0,\chi_1}(q^t), E_{\chi_1,\chi_0}(q^t) \mid t=1,2,4,8\}\cup\{B_1(q),B_2(q)\},\\
&\{E_{\chi_0,\chi_2}(q^t), E_{\chi_2,\chi_0}(q^t) \mid t=1,5\}\cup \{C_k(q)\mid k=1,2,3,4\}, \\
&\{E_{\chi_0,\chi_3}(q), E_{\chi_1,\chi_2}(q) ,E_{\chi_2,\chi_1}(q),E_{\chi_3,\chi_0}(q)\}\cup\{ D_k(q) \mid k=1,2,3,4 \}
\end{align*}
are bases for $M_2(\Gamma_0(40))$,  $M_2(\Gamma_0(40),\chi_1)$,  $M_2(\Gamma_0(40),\chi_2)$,  $M_2(\Gamma_0(40),\chi_3)$, respectively.
\end{theorem}

{\bf Proof.} The assertions follow from (2.2), Theorems 2.2 and 2.3.
\eop

\vspace{1mm}

We now give four theorems (Theorems 2.5--2.8) from which the theorems of Section 3 (Theorems 3.1--3.4) follow.

\begin{theorem} 
We have
\begin{align*}
&\begin{aligned}
\varphi^4 (q) =&  8L(q) -32L(q^4),
 \end{aligned}\\
 &\begin{aligned}
\varphi^2 (q)  \varphi^2 (q^2)=& 4 L(q) - 4 L(q^2) +8 L(q^4) - 32 L(q^{8}) ,
 \end{aligned}\\
 &\begin{aligned}
\varphi^2 (q)  \varphi^2 (q^5)=&  \frac{4}{3}L(q) - \frac{16}{3}L(q^4) + \frac{20}{3}L(q^5) - \frac{80}{3}L(q^{20}) 
    + \frac{8}{3}A_1(q),
 \end{aligned}\\
 &\begin{aligned}
\varphi^2 (q)  \varphi^2 (q^{10})=& \frac{2}{3}L(q) - \frac{2}{3}L(q^2) + \frac{4}{3}L(q^4) + \frac{10}{3}L(q^5) 
 - \frac{16}{3}L(q^8)  - \frac{10}{3} L(q^{10})  \\
 &  +\frac{20}{3}L(q^{20}) -\frac{80}{3}L(q^{40})+\frac{10}{3}  A_1(q) + \frac{8}{3}A_2(q) +4 A_3(q),
 \end{aligned}\\
 &\begin{aligned}
\varphi (q)\varphi (q^2)\varphi (q^5)\varphi (q^{10})=& L(q) - L(q^2) -2 L(q^4) -5L(q^5) +8 L(q^8) \\
&+5 L(q^{10})  +10 L(q^{20}) -40L(q^{40})  + A_1(q) + 2A_3(q),
 \end{aligned}\\
&\begin{aligned}
\varphi^2 (q^2)  \varphi^2 (q^5)=&  \frac{2}{3}L(q) - \frac{2}{3}L(q^2) + \frac{4}{3}L(q^4) + \frac{10}{3}L(q^5) 
 - \frac{16}{3}L(q^8)   - \frac{10}{3} L(q^{10})  \\
&+\frac{20}{3}L(q^{20}) -\frac{80}{3}L(q^{40}) -\frac{2}{3} A_1(q) + \frac{8}{3}A_2(q) -4 A_3(q).
 \end{aligned}
\end{align*}
\end{theorem}

{\bf Proof.} Let $(a_1,a_2,a_3,a_4)$ be any of the quaternary quadratic forms listed in  the first column of Table 2.1. By Theorem 2.1
we have $\varphi(q^{a_1})\varphi(q^{a_2})\varphi(q^{a_3})\varphi(q^{a_4}) \in M_2(\Gamma_0(40))$. 
By Theorem 2.4,  $\varphi(q^{a_1})\varphi(q^{a_2})\varphi(q^{a_3})\varphi(q^{a_4})$ must be a linear combination of $L(q)-tL(q^t)$ ($t=2,4,5,8,10,20,40$) 
and $A_k(q)$ ($k\in\{1,2,3\}$), namely 
\begin{align}
\varphi(q^{a_1})\varphi(q^{a_2})\varphi(q^{a_3})\varphi(q^{a_4}) =& x_1( L(q) -2L(q^2))+x_2 ( L(q) -4L(q^4))    \nonumber\\ 
&+ x_3( L(q) -5L(q^5)) + x_4 ( L(q) -8L(q^8))  \nonumber  \\
& + x_5 ( L(q) -10L(q^{10})) + x_6 ( L(q) -20L(q^{20}))   \\
& +x_7 ( L(q) -40L(q^{40})) +y_1 A_1(q) + y_2A_2(q) + y_3 A_3(q). \nonumber 
\end{align}
We only prove the last equation in the theorem as the others can be proven similarly. Let $(a_1,a_2,a_3,a_4)=(2,2,5,5)$.
Appealing to \cite[Theorem 3.13]{Kilford}, we find that the Sturm bound for the modular space $M_2(\Gamma_0(40))$ is $12$.
So, equating the coefficients of $q^{n}$ for $0\leq n\leq 12$ on both sides of (2.26),
we find a  system of linear equations with the unknowns $x_i$ ($1\leq i\leq 7$),  $y_1$, $y_2$ and $y_3$.
Using MAPLE we solve the system and find that
\begin{eqnarray*}
x_1=x_5=\frac{1}{3}, x_2=x_6=-\frac{1}{3}, x_3=y_1=-\frac{2}{3}, x_4=x_7=\frac{2}{3}, y_2=\frac{8}{3}, y_3=-4.
\end{eqnarray*}
Substituting these values back in (2.26), and with the obvious  simplifications, we find the asserted equation.
\eop

\vspace{2mm}
\begin{corollary} Let $n\in\nn$. We have
\begin{align*}
N(1,1,10,10;n)=N(2,2,5,5;n)  \mbox{ if } n\equiv 0 \pmdd 2.
\end{align*}
\end{corollary}

{\bf Proof.} From Theorem 2.5, we have 
\begin{align}
\varphi^2(q)\varphi^2(q^{10})-\varphi^2(q^2)\varphi^2(q^5)=4A_1(q)+8A_3(q).
\end{align}
It is clear from (1.2), (2.13) and (2.15) that
\begin{align}
a_1(n)= a_3(n)=0 \mbox{ if } n\equiv 0 \pmdd 2.
\end{align}
The assertion now follows from (1.1), (2.27) and (2.28). 
\eop
\vspace{2mm}

Similarly to Theorem 2.5, Theorems 2.6--2.8 follow from Theorems 2.1 and 2.4.

\begin{theorem}  
Let $\chi_0$ be the trivial character and $\chi_1$ be as in {\em(2.3)}. Then
\begin{align*}
&\begin{aligned}
\varphi^3 (q) \varphi (q^5)= & E_{\chi_0,\chi_1}(q) -2E_{\chi_0,\chi_1}(q^2) -4 E_{\chi_0,\chi_1}(q^4) + 5E_{\chi_1,\chi_0}(q)  \\
& + 10E_{\chi_1,\chi_0}(q^2) -20E_{\chi_1,\chi_0}(q^4),
\end{aligned}\\
&\begin{aligned}
\varphi^2 (q)\varphi (q^2)\varphi (q^{10})= &-\frac{1}{2}E_{\chi_0,\chi_1}(q) +\frac{1}{2} E_{\chi_0,\chi_1}(q^2) - E_{\chi_0,\chi_1}(q^4)  \\
 &-4 E_{\chi_0,\chi_1}(q^8)   +\frac{5}{2}E_{\chi_1,\chi_0}(q)+\frac{5}{2}E_{\chi_1,\chi_0}(q^2) \\
& + 5E_{\chi_1,\chi_0}(q^4) -20E_{\chi_1,\chi_0}(q^8)+ 2B_2(q),
\end{aligned}\\
&\begin{aligned}
\varphi(q)\varphi ^2(q^2)\varphi (q^5)= & \frac{1}{2}E_{\chi_0,\chi_1}(q) -\frac{1}{2} E_{\chi_0,\chi_1}(q^2) - E_{\chi_0,\chi_1}(q^4)  \\
 &  -4 E_{\chi_0,\chi_1}(q^8) +\frac{5}{2}E_{\chi_1,\chi_0}(q)+\frac{5}{2}E_{\chi_1,\chi_0}(q^2)\\
& - 5E_{\chi_1,\chi_0}(q^4)+20E_{\chi_1,\chi_0}(q^8)+5 B_1(q) -B_2(q),
\end{aligned}\\
&\begin{aligned}
\varphi (q) \varphi^3 (q^5)=& E_{\chi_0,\chi_1}(q) -2 E_{\chi_0,\chi_1}(q^2) -4 E_{\chi_0,\chi_1}(q^4) + E_{\chi_1,\chi_0}(q)  \\ 
&+ 2E_{\chi_1,\chi_0}(q^2) -4E_{\chi_1,\chi_0}(q^4),
\end{aligned}\\
&\begin{aligned}
\varphi(q)\varphi (q^5)\varphi^2 (q^{10})=&\frac{1}{2}E_{\chi_0,\chi_1}(q) -\frac{1}{2} E_{\chi_0,\chi_1}(q^2) - E_{\chi_0,\chi_1}(q^4) \\
 & -4E_{\chi_0,\chi_1}(q^8)+\frac{1}{2}E_{\chi_1,\chi_0}(q)+\frac{1}{2}E_{\chi_1,\chi_0}(q^2)\\
& - E_{\chi_1,\chi_0}(q^4)+4E_{\chi_1,\chi_0}(q^8)- B_1(q) +B_2(q),
\end{aligned}\\
&\begin{aligned}
\varphi(q^2)\varphi^2 (q^5)\varphi (q^{10})=&-\frac{1}{2}E_{\chi_0,\chi_1}(q) +\frac{1}{2} E_{\chi_0,\chi_1}(q^2) - E_{\chi_0,\chi_1}(q^4)\\
 &  -4E_{\chi_0,\chi_1}(q^8) +\frac{1}{2}E_{\chi_1,\chi_0}(q)+\frac{1}{2}E_{\chi_1,\chi_0}(q^2)\\
&+ E_{\chi_1,\chi_0}(q^4) -4E_{\chi_1,\chi_0}(q^8)-2 B_1(q).
\end{aligned}
\end{align*}
\end{theorem}

\begin{theorem} 
Let $\chi_0$ be the trivial character and  $\chi_2$ be as in {\em (2.3)}. Then
\begin{align*}
&\begin{aligned}
\varphi^3(q) \varphi (q^2) =& -2 E_{\chi_0,\chi_2}(q) + 8 E_{\chi_2,\chi_0}(q) ,
\end{aligned}\\
&\begin{aligned}
\varphi^2(q) \varphi (q^5) \varphi (q^{10})
=& \frac{2}{13}\big(2 E_{\chi_0,\chi_2}(q) - 15E_{\chi_0,\chi_2}(q^5) +8 E_{\chi_2,\chi_0}(q) + 60E_{\chi_2,\chi_0}(q^5)\big) \\
&+  \frac{8}{13}\big(6 C_1(q) -  4C_2(q) - 3 C_3(q) + 4 C_4(q)\big),
\end{aligned}\\
&\begin{aligned}
\varphi(q) \varphi^3 (q^2) =& -2 E_{\chi_0,\chi_2}(q) + 4 E_{\chi_2,\chi_0}(q) ,
&\end{aligned}\\
&\begin{aligned}
\varphi(q) \varphi (q^2) \varphi^2 (q^5)=& \frac{2}{13}\big(-3 E_{\chi_0,\chi_2}(q) -10E_{\chi_0,\chi_2}(q^5) +12E_{\chi_2,\chi_0}(q) \big)\\
& -\frac{80}{13}E_{\chi_2,\chi_0}(q^5) +  \frac{8}{13}\big(-2 C_2(q) -  5 C_3(q) +  C_4(q)\big),
\end{aligned}\\
 &\begin{aligned}
\varphi(q) \varphi (q^2) \varphi^2 (q^{10})=&\frac{2}{13}\big(-3 E_{\chi_0,\chi_2}(q) -10E_{\chi_0,\chi_2}(q^5) +6 E_{\chi_2,\chi_0}(q)\big)\\
&-\frac{40}{13}E_{\chi_2,\chi_0}(q^5) +  \frac{4}{13}\big(2 C_1(q) -  2 C_3(q) + 5 C_4(q)\big),
\end{aligned}\\
&\begin{aligned}
\varphi^2(q^2) \varphi (q^5) \varphi (q^{10})
=& \frac{2}{13}\big(2 E_{\chi_0,\chi_2}(q) -15E_{\chi_0,\chi_2}(q^5) +4 E_{\chi_2,\chi_0}(q) +30E_{\chi_2,\chi_0}(q^5)\big)   \\
& + \frac{4}{13}\big(-4 C_1(q) +  12 C_2(q) +  8 C_3(q) - 3 C_4(q)\big).
\end{aligned}
\end{align*}
\end{theorem}

\begin{theorem} 
Let $\chi_0$ be the trivial character and $\chi_1, \chi_2, \chi_3$ be as in {\em(2.2)}. Then
\begin{align*}
&\begin{aligned}
\varphi^3 (q)\varphi (q^{10})= & \frac{1}{7} \big(-E_{\chi_0,\chi_3}(q) - 5E_{\chi_1,\chi_2}(q) +4 E_{\chi_2,\chi_1}(q) 
+20E_{\chi_3,\chi_0}(q)\big)  \\
&+\frac{4}{7}\big(-3D_1(q)+15D_2(q) -15 D_3(q)+ 9 D_4(q)\big),
\end{aligned}\\
&\begin{aligned}
\varphi^2(q)\varphi (q^2)\varphi (q^5)=&\frac{1}{7}\big(- E_{\chi_0,\chi_3}(q) +5E_{\chi_1,\chi_2}(q) - 4 E_{\chi_2,\chi_1}(q) +20 E_{\chi_3,\chi_0}(q)\big)\\
& + \frac{8}{7}\big(-D_1(q) + 2D_4(q)\big),
\end{aligned}\\
&\begin{aligned}
\varphi(q)\varphi^2 (q^2)\varphi (q^{10})=& \frac{1}{7}\big(- E_{\chi_0,\chi_3}(q) -5E_{\chi_1,\chi_2}(q) +2 E_{\chi_2,\chi_1}(q) +10E_{\chi_3,\chi_0}(q)\big) \\
& +\frac{4}{7}\big(D_1(q)+ 5 D_2(q) + 5 D_3(q)+D_4(q)\big),
\end{aligned}\\
&\begin{aligned}
\varphi(q)\varphi^2 (q^5)\varphi (q^{10})=& \frac{1}{7} \big(-E_{\chi_0,\chi_3}(q) -E_{\chi_1,\chi_2}(q) +4 E_{\chi_2,\chi_1}(q) + 4 E_{\chi_3,\chi_0}(q)  \\
& +\frac{8}{7}\big(D_2(q) -D_3(q)+D_4(q)\big),
\end{aligned}\\
&\begin{aligned}
\varphi (q)\varphi^3 (q^{10})=& - \frac{1}{7} \big(E_{\chi_0,\chi_3}(q) + E_{\chi_1,\chi_2}(q) -2 E_{\chi_2,\chi_1}(q) - 2E_{\chi_3,\chi_0}(q) \big) \\
&+\frac{12}{7}\big(D_2(q) +D_3(q)+D_4(q)\big),
\end{aligned}\\
&\begin{aligned}
\varphi^3 (q^2) \varphi (q^5)=& \frac{1}{7}\big(- E_{\chi_0,\chi_3}(q) + 5E_{\chi_1,\chi_2}(q) - 2E_{\chi_2,\chi_1}(q) + 10E_{\chi_3,\chi_0}(q)
-12D_1(q)\big) ,
\end{aligned}\\
&\begin{aligned}
\varphi (q^2) \varphi^3 (q^5)=& \frac{1}{7}\big(- E_{\chi_0,\chi_3}(q) +E_{\chi_1,\chi_2}(q) -4 E_{\chi_2,\chi_1}(q) 
+4E_{\chi_3,\chi_0}(q) \big) \\
&-\frac{12}{7}\big(D_1(q)+D_2(q) + 3D_3(q) -D_4(q)\big),
\end{aligned}\\
&\begin{aligned}
 \varphi(q^2)\varphi (q^5)\varphi^2 (q^{10})=& \frac{1}{7}\big(- E_{\chi_0,\chi_3}(q) + E_{\chi_1,\chi_2}(q) - 2E_{\chi_2,\chi_1}(q)  + 2 E_{\chi_3,\chi_0}(q)\big)\\
&+\frac{4}{7}\big(-D_1(q)+D_2(q) - 3D_3(q)+D_4(q)\big).
\end{aligned}
\end{align*}
\end{theorem}

\section{Main results} 

\begin{theorem}  
Let $n \in \nn$.  We have 
\begin{align*}
&\begin{aligned}
N(1,1,5,5;n) = & \frac{4}{3}\sigma(n) - \frac{16}{3} \sigma (n/4) + \frac{20}{3} \sigma (n/5)  
-\frac{80}{3}\sigma(n/20)   + \frac{8}{3}a_1(n),
\end{aligned}\\
&\begin{aligned}
N(1,1,10,10;n) = & \frac{2}{3}\sigma(n) - \frac{2}{3} \sigma (n/2) + \frac{4}{3} \sigma (n/4) + \frac{10}{3} \sigma (n/5) - \frac{16}{3}\sigma(n/8)  \\ &- \frac{10}{3}\sigma(n/10) + \frac{20}{3}\sigma(n/20)  - \frac{80}{3}\sigma(n/40) + \frac{10}{3}a_1(n) \\ &+ \frac{8}{3}a_2(n) +4a_3(n),
\end{aligned}\\
&\begin{aligned}
N(1,2,5,10;n) = & \sigma(n) - \sigma (n/2) -2 \sigma (n/4) - 5 \sigma (n/5) + 8 \sigma(n/8)  \\ 
& + 5\sigma(n/10) + 10\sigma(n/20)  - 40\sigma(n/40) + a_1(n) + 2a_3(n),
\end{aligned}\\
&\begin{aligned}
N(2,2,5,5;n) = & \frac{2}{3}\sigma(n) - \frac{2}{3} \sigma (n/2) + \frac{4}{3} \sigma (n/4) + \frac{10}{3} \sigma (n/5) - \frac{16}{3}\sigma(n/8) \\ 
&- \frac{10}{3}\sigma(n/10)  +\frac{20}{3}\sigma(n/20)  - \frac{80}{3}\sigma(n/40) - \frac{2}{3}a_1(n) \\&+ \frac{8}{3}a_2(n) -4a_3(n).
\end{aligned}
\end{align*}
\end{theorem}

{\bf Proof.} The assertions follow from (1.1),  (2.4) and Theorem 2.5.
\eop

\begin{theorem}  
Let $n \in \nn$. Let $\sigma_{\chi_i, \chi_j}(n)$ be as in {\em(2.1)}  for $ i, j \in \{0,1\}$. We have  
\begin{align*}
&\begin{aligned}
N(1,1,1,5;n)= &\sigma_{\chi_0,\chi_1}(n) -2\sigma_{\chi_0,\chi_1}(n/2) -4 \sigma_{\chi_0,\chi_1}(n/4) \\
& + 5\sigma_{\chi_1,\chi_0}(n)+10\sigma_{\chi_1,\chi_0}(n/2)- 20\sigma_{\chi_1,\chi_0}(n/4),
\end{aligned}\\
&\begin{aligned}
N(1,1,2,10;n) = & -\frac{1}{2}\sigma_{\chi_0,\chi_1}(n) +\frac{1}{2} \sigma_{\chi_0,\chi_1}(n/2) - \sigma_{\chi_0,\chi_1}(n/4)
  -4 \sigma_{\chi_0,\chi_1}(n/8) \\
&+\frac{5}{2}\sigma_{\chi_1,\chi_0}(n)+\frac{5}{2}\sigma_{\chi_1,\chi_0}(n/2)+ 5\sigma_{\chi_1,\chi_0}(n/4)\\
&-20\sigma_{\chi_1,\chi_0}(n/8)  + 2b_2(n),
\end{aligned}\\
&\begin{aligned}
N(1,2,2,5;n)= &\frac{1}{2}\sigma_{\chi_0,\chi_1}(n) -\frac{1}{2} \sigma_{\chi_0,\chi_1}(n/2) - \sigma_{\chi_0,\chi_1}(n/4)
-4 \sigma_{\chi_0,\chi_1}(n/8)  \\
& +\frac{5}{2}\sigma_{\chi_1,\chi_0}(n)+\frac{5}{2}\sigma_{\chi_1,\chi_0}(n/2)- 5\sigma_{\chi_1,\chi_0}(n/4)
 +20\sigma_{\chi_1,\chi_0}(n/8) \\
 &+5 b_1(n) - b_2(n),
\end{aligned}\\
&\begin{aligned}
N(1,5,5,5;n)= &\sigma_{\chi_0,\chi_1}(n) -2\sigma_{\chi_0,\chi_1}(n/2) -4 \sigma_{\chi_0,\chi_1}(n/4) \\
& + \sigma_{\chi_1,\chi_0}(n)+2\sigma_{\chi_1,\chi_0}(n/2)- 4\sigma_{\chi_1,\chi_0}(n/4),
\end{aligned}\\
&\begin{aligned}
N(1,5,10,10;n)= &\frac{1}{2}\sigma_{\chi_0,\chi_1}(n) -\frac{1}{2} \sigma_{\chi_0,\chi_1}(n/2) - \sigma_{\chi_0,\chi_1}(n/4)
 -4\sigma_{\chi_0,\chi_1}(n/8) \\& +\frac{1}{2}\sigma_{\chi_1,\chi_0}(n)+\frac{1}{2}\sigma_{\chi_1,\chi_0}(n/2)- \sigma_{\chi_1,\chi_0}(n/4)
 +4\sigma_{\chi_1,\chi_0}(n/8) \\&- b_1(n) +b_2(n),
 \end {aligned}\\
&\begin{aligned}
N(2,5,5,10;n)=&-\frac{1}{2}\sigma_{\chi_0,\chi_1}(n) +\frac{1}{2} \sigma_{\chi_0,\chi_1}(n/2) - \sigma_{\chi_0,\chi_1}(n/4)
 -4\sigma_{\chi_0,\chi_1}(n/8)  \\
 &+\frac{1}{2}\sigma_{\chi_1,\chi_0}(n)+\frac{1}{2}\sigma_{\chi_1,\chi_0}(n/2)+ \sigma_{\chi_1,\chi_0}(n/4)\\
 & -4\sigma_{\chi_1,\chi_0}(n/8) -2 b_1(n).
\end {aligned}
\end{align*}
\end{theorem}

{\bf Proof.} The assertions follow from (1.1),  (2.5) and Theorem 2.6.
\eop

\begin{theorem} 
Let $n \in \nn$. Let $\sigma_{\chi_i, \chi_j}(n)$ be as in {\em(2.1)} for $ i, j \in \{0,2\}$. Then 
\begin{align*}
&\begin{aligned}
N(1, 1, 1,2;n) =&-2 \sigma_{\chi_0,\chi_2}(n) +8 \sigma_{\chi_2,\chi_0}(n),
\end{aligned}\\
&\begin{aligned}
N(1,1,5,10 ;n) = &  \frac{2}{13}\big(2 \sigma_{\chi_0,\chi_2}(n) -15\sigma_{\chi_0,\chi_2}(n/5) +8\sigma_{\chi_2,\chi_0}(n) + 60\sigma_{\chi_2,\chi_0}(n/5)\big) \\&
+  \frac{8}{13}\big(6 c_1(n) -  4 c_2(n) -  3 c_3(n) +  4 c_4(n)\big),
\end{aligned}\\
&\begin{aligned}
N(1, 2, 2,2;n) =&-2 \sigma_{\chi_0,\chi_2}(n) +4 \sigma_{\chi_2,\chi_0}(n),
\end{aligned}\\
&\begin{aligned}
N(1, 2, 5,5;n) =&-\frac{2}{13}\big(3 \sigma_{\chi_0,\chi_2}(n) +10\sigma_{\chi_0,\chi_2}(n/5) - 12 \sigma_{\chi_2,\chi_0}(n) +40\sigma_{\chi_2,\chi_0}(n/5) \big) \\
& + \frac{8}{13}\big(-2 c_2(n) -  5 c_3(n) +  c_4(n)\big),
\end{aligned}\\
&\begin{aligned}
N(1,2,10,10;n) = &\frac{2}{13}\big(-3 \sigma_{\chi_0,\chi_2}(n) -10\sigma_{\chi_0,\chi_2}(n/5) +6\sigma_{\chi_2,\chi_0}(n)  - 20\sigma_{\chi_2,\chi_0}(n/5)\big)\\ 
 &+  \frac{4}{13}\big(2 c_1(n) -  2 c_3(n) +  5c_4(n)\big),
\end{aligned}\\
&\begin{aligned}
N(2,2,5,10;n)= &\frac{2}{13}\big(2 \sigma_{\chi_0,\chi_2}(n) -15\sigma_{\chi_0,\chi_2}(n/5) +4 \sigma_{\chi_2,\chi_0}(n) +  30\sigma_{\chi_2,\chi_0}(n/5)\big)\\&
  + \frac{4}{13}\big(-4 c_1(n) + 12 c_2(n) +  8c_3(n) -  3 c_4(n)\big).
\end {aligned}
\end{align*}
\end{theorem}

{\bf Proof.} The assertions follow from (1.1),  (2.6) and Theorem 2.7.
\eop

\begin{theorem} 
Let $n \in \nn$. Let $\sigma_{\chi_i, \chi_j}(n)$ be as in {\em(2.1)}  for $ i, j \in \{0,1,2,3\}$. Then 
\begin{align*}
&\begin{aligned}
N(1,1,1,10;n)=&  \frac{1}{7}\big(- \sigma_{\chi_0,\chi_3}(n) -5\sigma_{\chi_1,\chi_2}(n) +4 \sigma_{\chi_2,\chi_1}(n) +20\sigma_{\chi_3,\chi_0}(n)\big)  
\\&+\frac{4}{7}\big(-3d_1(n)  +15d_2(n) -15 d_3(n)+ 9d_4(n)\big), 
\end{aligned} \\
&\begin{aligned}
N(1,1,2,5;n)= & \frac{1}{7}\big(- \sigma_{\chi_0,\chi_3}(n) +5\sigma_{\chi_1,\chi_2}(n) - 4 \sigma_{\chi_2,\chi_1}(n) + 20\sigma_{\chi_3,\chi_0}(n)\big)  
\\&+\frac{8}{7}\big(-d_1(n) + 2 d_4(n)\big), 
\end{aligned}\\
&\begin{aligned}
 N(1,2,2,10;n)= & \frac{1}{7}\big(- \sigma_{\chi_0,\chi_3}(n) -5\sigma_{\chi_1,\chi_2}(n) +2 \sigma_{\chi_2,\chi_1}(n) +10 \sigma_{\chi_3,\chi_0}(n) \big) 
\\&+\frac{4}{7}\big(d_1(n)  + 5d_2(n) + 5 d_3(n)+ d_4(n)\big), 
\end{aligned}\\
&\begin{aligned}
N(1,5,5,10;n)= &\frac{1}{7}\big(- \sigma_{\chi_0,\chi_3}(n) -\sigma_{\chi_1,\chi_2}(n) +4 \sigma_{\chi_2,\chi_1}(n) 
+4\sigma_{\chi_3,\chi_0}(n)\big)  \\
&+\frac{8}{7}\big(d_2(n)  -d_3(n)+d_4(n)\big), 
\end{aligned}\\
&\begin{aligned}
N(1,10,10,10;n) = &\frac{1}{7}\big(- \sigma_{\chi_0,\chi_3}(n) -\sigma_{\chi_1,\chi_2}(n) +2 \sigma_{\chi_2,\chi_1}(n) +2\sigma_{\chi_3,\chi_0}(n)\big)  \\
& +\frac{12}{7}\big(d_2(n) +d_3(n)+d_4(n)\big), 
\end{aligned}\\
&\begin{aligned}
N(2,2,2,5 ;n) = &  \frac{1}{7}\big(- \sigma_{\chi_0,\chi_3}(n) +5\sigma_{\chi_1,\chi_2}(n) -2 \sigma_{\chi_2,\chi_1}(n) +10\sigma_{\chi_3,\chi_0}(n) \big)\\
&-\frac{12}{7}d_1(n), 
\end{aligned}\\
&\begin{aligned}
N(2,5,5,5;n) =& \frac{1}{7}\big(- \sigma_{\chi_0,\chi_3}(n) +\sigma_{\chi_1,\chi_2}(n) -4 \sigma_{\chi_2,\chi_1}(n) 
+4\sigma_{\chi_3,\chi_0}(n)\big)  \\ 
&+\frac{12}{7}\big(-d_1(n) -d_2(n) -3d_3(n)+d_4(n)\big), 
\end{aligned}\\
&\begin{aligned}
N(2,5,10,10;n)= &\frac{1}{7}\big(- \sigma_{\chi_0,\chi_3}(n) +\sigma_{\chi_1,\chi_2}(n) -2 \sigma_{\chi_2,\chi_1}(n) +2\sigma_{\chi_3,\chi_0}(n)\big) \\ 
&+\frac{4}{7}\big(-d_1(n) +d_2(n) -3d_3(n)+d_4(n)\big). 
\end{aligned}
\end{align*}
\end{theorem}

{\bf Proof.} The assertions follow from (1.1),  (2.7), (2.8) and Theorem 2.8.
\eop

\section{Remarks} 

\begin{remark}
{\rm  
Replacing $q$ by $-q$ in $\varphi^3 (q) \varphi (q^5)$ in Theorem 2.6, we have
\begin{align}
\varphi^3 (-q) \varphi (-q^5) =& E_{\chi_0,\chi_1}(-q) -2E_{\chi_0,\chi_1}(q^2) -4 E_{\chi_0,\chi_1}(q^4)  \nonumber   \\
& + 5E_{\chi_1,\chi_0}(-q) + 10E_{\chi_1,\chi_0}(q^2) -20E_{\chi_1,\chi_0}(q^4).
\end{align}
Appealing to Theorem 2.3, we obtain 
\begin{align}
&E_{\chi_0 , \chi_1} (-q) = - E_{\chi_0 , \chi_1} (q) - 2E_{\chi_0 , \chi_1} (q^2) + 4 E_{\chi_0 , \chi_1} (q^4), \\
&E_{\chi_1 , \chi_0} (-q) = - E_{\chi_1 , \chi_0} (q) + 2E_{\chi_1 , \chi_0} (q^2) + 4 E_{\chi_1 , \chi_0} (q^4) .
\end{align} 
Substituting (4.2) and (4.3) in (4.1), we obtain
\begin{align}
&\varphi^3 (-q) \varphi (-q^5) =-E_{\chi_0,\chi_1}(q) -4E_{\chi_0,\chi_1}(q^2) - 5 E_{\chi_1,\chi_0}(q) + 20 E_{\chi_1,\chi_0}(q^2).
\end{align}
It can easily be seen that 
\begin{align}
&-E_{\chi_0,\chi_1}(q) -4E_{\chi_0,\chi_1}(q^2) = 1 + \sum_{n=1}^{\infty} \Big( \sum_{d\mid n} (-1)^d \dqu{5}{d} d \Big) q^n, \\
&- E_{\chi_1,\chi_0}(q) + 4 E_{\chi_1,\chi_0}(q^2) =  \sum_{n=1}^{\infty} \Big( \sum_{d\mid n} (-1)^d \dqu{5}{n/d} d \Big) q^n. 
\end{align} 
Now, appealing to (1.1) and (4.4)--(4.6), we obtain
\begin{align*}
&\sum_{n=0}^{\infty} N(1,1,1,5;n) (-q)^n
=\varphi^3 (-q) \varphi (-q^5) \nonumber \\
&\hspace{20mm}= 1 + \sum_{n=1}^{\infty} \Big( \sum_{d\mid n} (-1)^d \dqu{5}{d} d \Big) q^n    
+ 5   \sum_{n=1}^{\infty} \Big( \sum_{d\mid n} (-1)^d \dqu{5}{n/d} d \Big) q^n,
\end{align*}
from which we deduce 
\begin{align*}
N(1,1,1,5;n)=  \sum_{d\mid n} (-1)^{n+d} \dqu{5}{d} d  + 5  \sum_{d\mid n} (-1)^{n+d} \dqu{5}{n/d} d,
\end{align*} 
which agrees with known results, see for example \cite[Theorem 5.1]{alaca2007quaternary}. 
Similarly, one can show that our formula for $N(1,5,5,5;n)$ given in Theorem 3.2 agrees with the result in
\cite[Theorem 6.1]{alaca2007quaternary}.
}
\end{remark}

\begin{remark}
{\rm
Appealing to Lemma 2.1 and Theorem 2.3, we obtain the following identities:
\begin{align*}
&L(q)-4L(q^4)=\frac{1}{8}\,   \frac{ \eta^{20} (2z) }{\eta^8 (z) \eta^8(4z)} ,\\  
& E_{\chi_0,\chi_1}(q) =-\frac{1}{5} \,  \frac{ \eta^5 (z) }{\eta (5z)} ,\\
& E_{\chi_1,\chi_0}(q) =  \frac{ \eta^5 (5z) }{\eta (z)}  ,\\
& E_{\chi_0,\chi_2}(q) =- \frac{1}{2} \,  \frac{ \eta^2(z) \eta(2z) \eta^3 (4z)}{\eta^2 (8z)}  ,\\
& E_{\chi_2,\chi_0}(q) = \frac{ \eta^3(2z) \eta(4z) \eta^2 (8z)}{\eta^2 (z)}  ,\\
& E_{\chi_0,\chi_1}(q) +4E_{\chi_0,\chi_1}(q^2)=- \frac{ \eta (z) \eta^2 (2z) \eta^3 (5z)}{\eta^2 (10z)}  ,\\
& E_{\chi_1,\chi_0}(q) +E_{\chi_1,\chi_0}(q^2)= \frac{ \eta^3 (2z) \eta^2 (5z) \eta(10z)}{\eta^2(z)}  ,\\
& E_{\chi_1,\chi_0}(q) - 4E_{\chi_1,\chi_0}(q^2)=  \frac{ \eta^3 (z) \eta (5z) \eta ^2 (10z)}{\eta^2 (2z)} , \\
& E_{\chi_0,\chi_2}(q) - 2E_{\chi_2,\chi_0}(q)= -\frac{1}{2}\, \frac{ \eta^{13} (4z) }{\eta^2 (z) \eta(2z) \eta^6 (8z)},\\
& E_{\chi_0,\chi_2}(q) - 4E_{\chi_2,\chi_0}(q)= -\frac{1}{2}\, \frac{ \eta^{13} (2z) }{\eta^6 (z) \eta(4z) \eta^2 (8z)}, \\
& E_{\chi_0,\chi_1}(q) -2E_{\chi_0,\chi_1}(q^2)-4E_{\chi_0,\chi_1}(q^4)= \frac{ \eta^5 (2z) \eta^7 (10z)}{\eta(z)\eta(4z) \eta^3 (5z) \eta^3 (20z)}  ,\\
& E_{\chi_1,\chi_0}(q) +2E_{\chi_1,\chi_0}(q^2)-4E_{\chi_1,\chi_0}(q^4)= \frac{ \eta^7 (2z) \eta^5 (10z)}{\eta^3(z)\eta^3(4z) \eta (5z) \eta (20z)}.
\end{align*}
}
\end{remark}

\begin{remark}
{\rm 
Set $a:=\varphi(q)$,  $b:=\varphi(q^2)$,  $c:=\varphi(q^5)$ and $d:=\varphi(q^{10})$.
We obtain the following identities from Theorem 2.8:
\begin{align*}
&ad(-a^2-b^2+5c^2-5d^2)+bc(5a^2-8b^2-5c^2+10d^2)=12D_1(q),\\
&ad(2a^2-b^2-4c^2+d^2)+bc(-a^2+b^2-5c^2+7d^2)=24D_2(q),\\
&ad(-a^2+2b^2-7c^2+10d^2)+bc(-a^2+4b^2+c^2-8d^2)=48D_3(q),\\
&ad(a^2-8b^2-5c^2+20d^2)+bc(7a^2-10b^2+5c^2-10d^2)=48D_4(q).
\end{align*}
}
\end{remark}

\begin{remark}
{\rm 
It would be interesting to determine  general formulas for the number of representations of a positive integer $n$ by the quaternary quadratic forms with coefficients in $\{1, p, q,  pq\}$, 
where $p$ and $q$ are distinct prime numbers. 
The case when $p=2$ and $q=7$ is treated in \cite{Ayse-Jamilah}. \\
}
\end{remark}

\noindent{\bf Acknowledgements.} 
The research of the first author was supported by a Discovery Grant
from the Natural Sciences and Engineering Research Council of Canada (RGPIN-418029-2013). 


\vspace{3mm}

\noindent
Ay\c{s}e Alaca and Mada Altiary \\
School of Mathematics and Statistics \\ 
Carleton University \\ 
Ottawa, ON K1S 5B6, Canada \\

\noindent
AyseAlaca@cunet.carleton.ca \\
MadaAltiary@cmail.carleton.ca

\end{document}